\pdfoutput=1

\documentclass[11pt,a5paper,twoside,fleqn,leqno,color,final,lm,md]{pre} 

\usepackage[applemac]{inputenc}
\usepackage{amsthm,math,pinlabel,accents}
\usepackage[color]{showkeys}
\usepackage{multicol}
\usepackage{microtype}

\let\a\alpha
\let\d\delta
\let\e\epsilon
\let\k=\kappa
\let\l\lambda
\let\r\rho
\let\s\sigma
\let\t\theta
\let\p\varphi

\def\W  {O}          
\def\Oma{\Om_\al}    
\def\H  {\bar H}     
\def\Nh {\smash{\hat N}}
\def\lip{\Ls}
\def\omt{\tilde\om}
\def\Fsj{\Fs^{\j}}
\def\dth{\upd\t}

\def\j{\ifmt\!\fi j}

\def\abs#1{\absabs#1::/}
\def\absabs#1:#2:#3/{\n{#1}\def\next{#2}\ifx\next\empty\else_{#2}\fi}
\catcode`\!=\active \let!\Bar

\def\eval#1\at{\rbar{#1}_}
\def\po{\putdelims\lbrace\rbrace\doublearg}

\def\ham.{hamiltonian}
\def\suff.{sufficiently}
\def\nbhd.{neighbourhood}
\def\wrt.{with respect to}
\def\resp.{respectively}
\def\co.{coordinate}
\def\trans.{trans\-for\-ma\-tion}
\def\Po.{Poincar\'e}
\def\kam.{{\scshape Kam}}
\def\ra.{real analytic}
\def\vf.{vector field}

\theoremstyle{plain}
\newtheorem{prp}{Proposition}
\newtheorem{lem}{Lemma}
\newtheorem{alm}{Lemma}[section]

\begin{document}

\title  {A Lecture\\
         on the\\
         Classical KAM Theorem}
\author {Jürgen Pöschel}
\date   {Version 1.3 -- 2009\\
         Small corrections\\ and updated references}

\address{Institut für Analysis, Dynamik und Optimierung\\
         Universität Stuttgart, Pfaffenwaldring 57, D-70569 Stuttgart\\
         poschel@mathematik.uni-stuttgart.de}

\maketitle

\section{The Classical KAM-Theorem}

\subdiv{a}

The purpose of this lecture is to describe the \kam. theorem in
its most basic form and to give a complete and detailed proof.
This proof essentially follows the traditional lines laid out by
the inventors of this theory, {\sl K}olmogorov, {\sl A}rnold and
{\sl M}oser (whence the acronym `\kam.'), and the emphasis is more
on the underlying ideas than on the sharpness of the arguments.
After all, \kam. theory is not only a collection of specific
theorems, but rather a methodology, a collection of ideas of how
to approach certain problems in perturbation theory connected
with ‘small divisors’.

\subdiv{b}

The classical \kam. theorem is concerned with the stability of
motions in \ham. systems, that are small perturbations of
integrable \ham. systems. These integrable systems are
characterized by the existence of action angle coordinates such
that the \ham. depends on the action variable alone -- see
\cite{AMM,MZ} for details. Thus we are going to consider
hamiltonians of the form
\[
  H(p,q) = h(p) + f_\ep(p,q), \qquad
  f_\ep(p,q) = \e f_*(p,q,\e)
\]
for small~$\e$, where $p=(p_1,\dots,p_n)$ are the action
variables varying over some domain $D\subset\R^n$, while
$q=(q_1,\dots,q_n)$ are the conjugate angular variables, whose
domain is the usual $n$-torus~$\T^n$ obtained from $\R^n$ by
identifying points whose components differ by integer multiples
of~$2\pi$. Thus, $f_\ep$ has period~$2\pi$ in each component
of~$q$. Moreover, all our \ham.s are assumed to be \emph{real
analytic} in all arguments.

The equations of motion are, as usual,
\[
  \dot p = -H_q(p,q), \qquad \dot q = H_p(p,q)
\]
in standard vector notation, where the dot indicates
differentiation with respect to the time~$t$, and the subscripts
indicate partial derivatives. The underlying phase space is
$D\x\T^n \subset \R^n\x\T^n$ with the standard symplectic
structure
\[
  \ups = \sum_{1\le j\le n} dp_j\wedge dq_j.
\]
The \ham. \vf.~$X_H$ associated with the equations
of motions then satisfies $\ups(X_H,\cd) = -dH$.

We assume that the number~$n$ of degrees of freedom is at
least~2, since one degree of freedom systems are always
integrable.

\subdiv{c}

For $\e=0$ the system is governed by the unperturbed, integrable
\ham.~$h$, and the equations of motion reduce to
\[
  \dot p = 0, \qquad \dot q = \om
\]
with
\[
  \om = h_p(p).
\]
They are easily integrated -- hence the name
\emph{integrable system} -- and their general solution is
\[
  p(t) = p_0, \qquad q(t) = q_0+\om(p_0)t.
\]
Hence, every solution curve is a straight line, which, due to
the identification of the $q$-coordinates modulo~$2\pi$, is
winding around the \emph{invariant torus} 
\[
  \Ts_{p_0} = \{p_0\}\x\T^n
\]
with constant {frequencies}, or winding numbers,
$\om(p_0)=(\om_1(p_0),\dots,\om_n(p_0))$. Such tori with linear
flow are also called \emph{Kronecker tori}.

In addition these tori are \emph{Lagrangian}. That is, the
restriction of the symplectic form~$\ups$ to their tangent space
vanishes, and their dimension is maximal \wrt. this property.

Thus, the whole phase space is \emph{foliated} into an
$n$-parameter family of invariant Lagrangian tori $\Ts_{p_0}$,
on which the flow is linear with constant {frequencies} $\om(p_0)$.  
-- This is the geometric picture of an integrable \ham. system.

It should be kept in mind that due to the introduction of
action angle \co.s these solutions are related to ‘real world
solutions’ by some \co. transformation, which is periodic in
$q_1,\dots,q_n$. Expanding such a transformation into Fourier
series and inserting the linear solutions obtained above, the
‘real world solutions’ are represented by series of the form
\[
  \Xi(p_0,q_0+\om(p_0)t)
  = \sum_{k\in\Z^n} a_k(p_0)\,e^{\pair k{q_0}+t\pair k{\om(p_0)}},
  \qquad a_k \in\R^{2n},
\]
where $\pair\cd\cd$ denotes the usual scalar product.
Thus, every solution is now \emph{quasi-periodic} in~$t$: its
frequency spectrum in general does not consist of integer
multiples of a single frequency -- as is the case with periodic
solutions --, but rather of integer combinations of a finite
number of different frequencies. In essence, the ‘real world
solutions’ are superpositions of $n$~oscillations, each 
with its own frequency. Moreover, these quasi-periodic solutions
occur in families, depending on the parameter~$q_0$, which
together fill an invariant embedded $n$-torus.

Let us return to action angle \co.s. We observe that the
topological nature of the flow on each Kronecker torus
crucially depends on  the arithmetical properties
of its frequencies~$\om$. There are essentially two cases.

1.~--~The frequencies $\om$ are \emph{nonresonant}, or
\emph{rationally independent}:
\[
  \pair k\om  \ne 0 \qtext{for all} 0\ne k\in\Z^n.
\]
Then, on this torus, each orbit is dense, the flow 
is ergodic, and the torus itself is minimal.

2.~--~The frequencies $\om$ are \emph{resonant}, or
\emph{rationally dependent}: that is, there exist integer
relations 
\[
  \pair k\om  = 0 \qtext{for some} 0\ne k\in\Z^n.
\]
The prototype is $\om=(\om_1,\dots,\om_{n-m},0\dots,0)$, with
$1\le m \le n-1$ trailing zeroes and nonresonant
$(\om_1,\dots,\om_{n-m})$. In this case the torus decomposes into
an $m$-parameter family of invariant $n-m$-tori.
Each orbit is dense on such a lower dimensional
torus, but not in~$\T^n$.
 
A special case arises when there exist $m=n-1$ independent
resonant relations. Then each frequency $\om_1,\dots,\om_n$ is an
integer multiple of a fixed non-zero frequency~$\om_*$, and the
whole torus is filled by \emph{periodic} orbits with one and
the same period~$2\pi/\om_*$.

In an integrable system the frequencies on the tori may or may
not vary with the torus, depending on the nature of the
\emph{frequency map}
\[
  h_p\maps D \to \Om, \quad p\mapsto \om(p)=h_p(p),
\]
where $\Om\subset\R^n$. We now make the assumption that this
system is \emph{nondegenerate} in the sense that
\[
  \det h_{pp} = \det {\del \om\over\del p} \ne 0
\]
on~$D$. Then $h_p$ is an open map, even a local diffeomorphism
between $D$ and some open \emph{frequency domain}
$\Om\subset\R^n$, and ‘the frequencies~$\om$ effectively depend
on the amplitudes~$p$’, as a physicist would say. It follows
that nonresonant tori and resonant tori of all types each form
dense subsets in phase space. Indeed, the resonant ones sit
among the nonresonant ones like the rational numbers among the
irrational numbers.

This ‘frequency-amplitude-modulation’ is a genuinely nonlinear
phenomenon. By contrast, in a linear system the frequencies are
the same all over the phase space. As we will see, this is
essential for the stability results of the \kam. theory. As
it is said, ‘the nonlinearities have a stabilizing effect’.

\subdiv{d}

Now we consider the perturbed \ham.. The objective is to prove
the persistence of invariant tori for small~$\e\ne0$.

The first result in this direction goes back to \Po. and is of a
negative nature. He observed that the \emph{resonant} tori are in
general \emph{destroyed} by an arbitrarily small perturbation.
In particular, out of a torus with an $n-1$-parameter family of
periodic orbits, usually only \emph{finitely} many periodic orbits
survive a perturbation, while the others disintegrate and give
way to \emph{chaotic behavior}. -- So in a nondegenerate system
a \emph{dense} set of tori is usually destroyed. This, in
particular, implies that a generic \ham. system is \emph{not
integrable} \cite{Ge,MM,Ze73}.

Incidentally, it would not help to drop the nondegeneracy
assumption to avoid resonant tori. If $h$ is too degenerate, the
motion may even become ergodic on each energy surface, thus
destroying all tori~\cite{Kat}.

A dense set of tori being destroyed there seems to be no hope
for other tori to survive. Indeed, until the fifties it was a
common belief that arbitrarily small perturbations can turn an
integrable system into an ergodic one on each energy surface.
In the twenties there even appeared an -- erroneous -- proof of
this ‘ergodic hypothesis’ by Fermi.

But in {1954} Kolmogorov observed that the converse is
true -- the \emph{majority} of tori survives. He proved the
persistence of those Kronecker systems, whose frequencies~$\om$
are not only nonresonant, but are \emph{strongly nonresonant} in
the sense that there exist constants $\a>0$ and $\tau>0$
such that
\[
  \abs{\pair k\om} \ge {\a\over\abs k^\tau}
  \qquad\text{for all $0\ne k\in\Z^n$},           
  \eqlabel{Dio}
\]
where $\abs k = \abs{k_1}+\dots+\abs{k_n}$. Such a condition is
called a \emph{diophantine} or \emph{small divisor condition}.

The existence of such frequencies is easy to see. Let
$\Dl_\a^\tau$ denote the set of all $\om\in\R^n$ satisfying
these infinitely many conditions with fixed $\a$ and~$\tau$. Then
$\Dl_\a^\tau$ is the complement of the open dense set
\[
  R_\a^\tau = \bigcup_{0\ne k\in\Z^n} R_{\a,k}^\tau,
\]
where
\[
  R_{\a,k}^\tau =  
  \set{\om\in\R^n: \abs{\pair k\om}<{\a/\abs k^\tau}}.
\]
Obviously, for any bounded domain $\Om\subset\R^n$, we have the
Lebesgue measure estimate
$m(R_{\a,k}^\tau\cap\Om) = O(\a/\abs k^{\tau+1})$, 
and thus
\[
  m(R_a^\tau\cap\Om)  \le \sum_{k\ne0} m(R_{\a,k}^\tau\cap\Om)
  =\bigo\a,
\]
provided that $\tau>n-1$. Hence, $ R^\tau = \bigcap_{\a>0}
R_\a^\tau $ is a set of measure zero, and its complement
\[
  \Dl^\tau = \bigcup_{\a>0} \Dl_\a^\tau
\]
is a set of \emph{full measure} in $\R^n$, for any $\tau>n-1$.
In other words, \emph{almost every} $\om$ in~$\R^n$
belongs to $\Dl^\tau$, $\tau>n-1$, which is the set of
all $\om$ in~$\R^n$ satisfying \eqref{Dio} for some $\a>0$
while $\tau$ is fixed..

As an aside we remark that $ \Dl^\tau = \varnothing$ for
$\tau<n-1$, because for every nonresonant~$\om$,
\[
  \min_{0\ne\abs{k}_\infty\le K} \abs{\pair k\om} 
  \le {\abs\om\over K^{n-1}}           
  \eqlabel{Dir}
\]
by Dirichlet's pigeon hole argument. And for $\tau=n-1$, the set
$\Dl^{n-1}$ has measure zero, but Hausdorff
dimension~$n$ -- see \cite{RuDS} for references. So there are
continuum many diophantine frequencies to the exponent~$n-1$,
although they form a set of measure zero.
-- 
But here we will fix $ \tau > n-1 $ and drop it from the
notation, letting $\Dl_\a = \Dl_\a^\tau$.

\subdiv{e}

But although almost all frequencies are strongly nonresonant for
any fixed $\tau>n-1$, it is \emph{not true} that almost all tori
survive a given perturbation~$f_\ep$,
no matter how small~$\e$. The reason is that the parameter~$\a$
in the nonresonance condition limits the size of the
perturbation through the condition
\[
  \e \ll \a^2.
\]
Conversely, under a given small perturbation of size~$\e$, only
those Kronecker tori with frequencies~$\om$ in $\Dl_\a$ with
\[
  \a \gg \sqrt\e,
\]
do survive. Thus, we can not allow $\a$ to vary, but have to
\emph{fix} it in advance.

To state the \kam. theorem, we therefore single out the
subsets
\[
  \Oma \subset \Om, \qquad \a>0,
\]
whose frequencies belong to $\Dl_\a$ and also have at
least distance~$\a$ to the boundary of~$\Om$. These, like
$\Dl_\a$, are \emph{Cantor sets}: they are closed, perfect and
\emph{nowhere dense}, hence of first Baire category. But they
also have \emph{large Lebesgue measure}: 
\[
  m(\Om-\Oma) = \bigo\a,
\]
provided the boundary of~$\Om$ is piecewise smooth, or at least
of dimension $n-1$  so that the measure of a boundary layer of
size~$\a$ is $\bigo\a$. 

The main theorem of Kolmogorov, Arnold and Moser can now be
stated as follows.

\begin{proclaim}[The Classical \kam. Theorem (\cite{A63,Kol,Mo63})]
Suppose the integrable \ham.~$h$ is nondegenerate, the
frequency map~$h_p$ is a diffeomorphism $D\to\Om$, and $H=h+f_\ep$
is real analytic on $\bar D\x\T^n$. Then there exists a
constant $\d>0$ such that for
\[
  \abs\e < \d\a^2
\]
all Kronecker tori $(\T^n,\om)$ of the unperturbed system with
$\om\in\Oma$ persist as Lagrangian tori, being only slightly
deformed. Moreover, they depend in a Lipschitz continuous way
on~$\om$ and fill the phase $D\x\T^n$ up to a set of
measure~$\bigo{\a}$. 
\end{proclaim}

Here, ‘real analytic on $\bar D\x\T^n$’ means that the
analyticity extends to a uniform \nbhd. of~$D$.

Is is an immediate and important consequence of the \kam. theorem
that small perturbations of nondegenerate \ham.s are \emph{not
ergodic}, as the Kronecker tori form an invariant
set, which is neither of full nor of zero measure. Thus the
ergodic hypothesis of the twenties was wrong.

It has to be stated again, however, that this invariant set,
although of large measure, is a \emph{Cantor set} and thus has no
interior points. It is therefore impossible to tell with finite
precision whether a given initial position falls onto an
invariant torus or into a gap between such tori. From a physical
point of view the \kam. theorem rather makes  a
\emph{probabilistic} statement: with probability $1-\bigo\a$ a
randomly chosen orbit lies on an invariant torus and is thus
perpetually stable.

\subdiv{f}

We conclude with some remarks about the necessity of the
assumptions of the \kam. theorem.

First, neither the perturbation nor the integrable \ham. need to
be real analytic. It suffices that they are differentiable of
class~$C^l$ with
\[
  l > 2\tau+2 > 2n
\]
to prove the persistence of individual tori \cite{Mo69,Po1,Sal}.
For their Lipschitz dependence some more regularity is
required~\cite{Po2}.

The nondegeneracy condition may also be relaxed. It is not
necessary that the frequency map is open. Roughly speaking, it
suffices that the intersection of its range with any hyperplane
in $\R^n$ has measure zero. For example, if it happens that
\[
  h_p(p) = (\om_1(p_1),\dots,\om_n(p_1))
\]
is a function of $p_1$ only  (and thus completely degenerate), it
suffices to require that
\[
  \det\pas{\del^j\om_i\over\del p_1^j}_{1\le i,j\le n} \ne 0.
\]
For a more general statement see~\cite{RuDV}, and \cite{Ru98,Sev,XYQ} for
proofs.

Finally, the \ham. nature of the equations is almost
indispensable. Analogous result are true for reversible systems
\cite{MoR,Po2}. But in any event the system has to be
\emph{conservative}. Any kind of dissipation immediately destroys
the Cantor family of tori, although isolated ones may persist as
attractors.

\section{The KAM Theorem with Parameters}

\subdiv{a}

Instead of proving the classical \kam. theorem directly, we are
going to deduce it from another \kam. theorem, which is concerned
with perturbations of a family of linear \ham.s. This is
accomplished by introducing the frequencies of the Kronecker
tori as independent parameters. This approach was first taken
in~\cite{Mo67}.
 
To this end we write $p=p_0+I$ and expand $h$ around~$p_0$
so that
\[
  h(p) 
  = h(p_0) + \pair*{h_p(p_0)}I 
  + \int_0^1 (1-t)\pair*{h_{pp}(p_t)I}I\dt,
\]
where $p_t=p_0+tI$. By assumption,
the frequency map is a \emph{diffeomorphism}
\[
  h_p \maps D \to \Om, \quad p_0 \mapsto \om=h_p(p_0).
\]
Hence, instead of $p_0\in D$ we may introduce the frequencies
$\om\in\Om$ as independent parameters, determining $p_0$
uniquely.
--
Incidentally, the inverse map is given as
\[
  g_\om\maps \Om \to D, \quad \om \mapsto p_0=g_\om(\om),
\]
where $g$ is the \emph{Legendre transform} of~$h$, defined by
$g(\om) = \sup_p (\pair p\om-h(p))$.
See \cite{AMM} for more details on Legendre transforms.

Thus we write
\[
  h(p)=e(\om)+\pair\om I + P_h(I,\om)
\]
with
\[
  P_h = \int_0^1 (1-t)\pair*{h_{pp}(p_t)I}I\dt,
\]
and
\[
  f_\ep(p,q)
  = P_{\e}(I,q,\om)
  = \e f_*(p_0+I,q,\e) .
\]
Writing also $\t$ instead of~$q$ for the angular variables, we
obtain the family of \ham.s $H=N+P$ with
\[
  N = e(\om) +\pair\om I, \quad
  P = P_h(I,\om) + P_{\e}(I,\t,\om).
\]
They are real analytic in the \co.s $(I,\t)$ in $B\x\T^n$, $B$
some \suff. small ball around the origin in~$\R^n$, as well as
in the parameters~$\om$ in some \emph{uniform} \nbhd. of~$\Om$. 

This family is our new starting point. For $P=0$ it
reduces to the \emph{normal form} 
$N = e(\om)+\pair\om I$. There is an invariant Kronecker torus
\[
  \Ts_\om = \{0\}\x\T^n
\]
with constant \vf. 
\[
  X_N = \sum_{1\le j\le n} \om_j\,{\del\hfil\over\del\t_j}
\] 
for each $\om\in\Om$, and all these tori are given by the
family 
\[
  \Phi_0\maps \T^n\x\Om \to B\x\T^n, \quad
  (\t,\om) \mapsto (0,\t)
\]
of trivial embeddings of $\T^n$ over~$\Om$ into phase space.
Moreover, each such torus clearly is Lagrangian.
Our aim is to prove the persistence of a `subfamily' of such
Lagrangian Kronecker tori  under \suff. small perturbations
$P\ne0$ over the Cantor set $\Oma\subset\Om$ of
frequency parameters~$\om$. -- Thus, instead of proving the
existence of a {Cantor family} of invariant tori in \emph{one
\ham. system}, we first prove the existence of \emph{one invariant
torus} within a Cantor family of \ham. systems.

This change of perspective has several advantages. --- The
unperturbed \ham. is as simple as possible, namely linear. This
simplifies the \kam. proof. --- The frequencies are separated from
the actions. This makes their r\^ole more transparent. For
example, the Lipschitz dependence of the tori on~$\om$ is easily
established. --- Generalizations such as weaker nondegeneracy
conditions and extension to infinite dimensional systems are
easier. Also, this approach lends itself to applications in
bifurcation theory, where systems naturally depend on parameters.

\subdiv{b}

To state the basic result quantitatively we need to introduce a
few notations. Let 
\[
  D_{r,s} = \{I:\abs{I}<r\}\x\{\t:\abs{\Im\t}<s\} 
  \subset \C^n\x\C^n 
\] 
and 
\[
  \W_h = \{\om:\abs{\om-\Oma}<h\} \subset \C^n
\]
denote complex \nbhd.s of the torus $\{0\}\x\T^n$ and~$\Oma$,
\resp., where $\abs\cd$ stands for the sup-norm of real
vectors. The sup-norm of functions on $D_{r,s}\x\W_h$ is denoted
by $\abs{\cd:r,s,h}$.

We will also consider the Lipschitz constants of mappings
\wrt.~$\om$. We define  
\[ 
  \abs{\p:\lip} = \sup_{\ups\ne\om} 
    {\abs{\p(\ups)-\p(\om)}\over\abs{\ups-\om}},
\]
where the \iflm underlying\fi domain will be clear from the context or
indicated by a subscript.

\begin{proclaim}[Theorem A] 
Let $H=N+P$. Suppose $P$ is real analytic on $D_{r,s}\x\W_h$ with
\[
  \abs{P:r,s,h} \le \gm\a rs^\nu, \qquad \a s^\nu\le h,
\]
where $\nu=\tau+1$ and $\gm$ is a small constant depending only
on $n$ and~$\tau$. Suppose also that
$r,s,h\le1$. Then there exists a Lipschitz continuous map
$
\p\maps \Oma \to \Om
$
close to the identity and a Lipschitz continuous family of \ra.
torus embeddings
$
\Phi\maps \T^n\x\Oma \to B\x\T^n
$
close to $\Phi_0$, such that for each $\om\in\Oma$ the embedded
tori are Lagrangian and
\[
  \rbar{X_H}_{\p(\om)} \comp \Phi = D\Phi\cd X_N.
\]
Moreover, $\Phi$ is \ra. on
$T_*=\{\t:\abs{\Im\t}<s/2\}\,$ for each~$\om$, and
\begin{align*}
  \abs{W(\Phi-\Phi_0)}, \a s^\nu\abs{W(\Phi-\Phi_0):\lip} 
  &\le {c\over\a rs^\nu} \abs{P:r,s,h},
  \\
  \abs{\p-\id}, \a s^\nu\abs{\p-\id:\lip} 
  &\le {c\over r}\abs{P:r,s,h},
\end{align*}
uniformly on $T_*\x\Oma$ and $\Oma$, \resp., where $c$ is a
constant depending only on $n$ and~$\tau$, and
$W=\diag(r\inv\Id,s\inv\Id)$.
\end{proclaim}

By slight abuse of notation we wrote $X_N$ also for the
analogous constant \vf. on~$\T^n$ alone. --
The theorem states that for each $\om$ in~$\Oma$ there is an
embedded invariant Kronecker torus $\Ts_\om=\Phi(\T^n,\om)$ with
frequencies~$\om$ for the \ham. \vf.~$X_H$ at the slightly
\emph{shifted} parameter value
$\omt = \p(\om)$.
Conversely, for each $\omt$ in the slightly deformed
Cantor set
\[
  \tilde\Om_\al = \p(\Oma) \subset \Om,
\]
the \vf. $\rbar1{X_H}_{\omt}$ admits an invariant Kronecker torus
$\Ts_\om$ with frequencies $\om=\p\inv(\omt)$. Each
such torus is Lagrangian and close to the
corresponding unperturbed torus.

\subdiv{c}

The Lipschitz estimates allow us to control the measure of
$\tilde\Om_\al$ and its complement. To this end we first extend $\p$ to
a lipeomorphism $!\p$ of~$\Om$.

\begin{prp}
The map $\p$ can be extended to a lipeomorphism
$!\p\maps\Om\to\Om$ with \vglue-1.5pc
\[
  \abs{!\p-\id:\lip,\Om} 
  \le \max\pas{\abs{\p-\id:\lip,\Oma},\a\inv\abs{\p-\id:\Oma}}.
\]
\end{prp}

\begin{proof}
Let $\psi$ be any \co. function of $\p-\id$ defined on~$\Oma$.
Define $\psi$ on $\Om^c=\R^n-\Om$ to be zero. Then
\[
  \abs{\psi:\lip,\Oma\cup\Om^c} \le 
  \max\pas{\abs{\p-\id:\lip,\Oma},\a\inv\abs{\p-\id:\Oma}} < 1.
\]
According to Appendix~B we can extend $\psi$ to a function
$!\psi$ on all of~$\R^n$ preserving its Lipschitz constant .
Doing this with every coordinate of~$\p$ we obtain an extension
$!\p$ of $\p$ such that $!\p=\id$
on~$\Om^c$ and
\[
  \abs{!\p-\id:\lip,\R^n} = \abs{\psi:\lip,\Oma\cup\Om^c} < 1.
\]
Hence $!\p$ is a lipeomorphism
on~$\R^n$. Since it is the identity outside of~$\Om$, it is
also a lipeomorphism of~$\Om$, extending~$\p$.\qed
\end{proof}

\begin{prp}
For $\tilde\Om_\al=\p(\Oma)$ the estimate
\[
  m(\Om-\tilde\Om_\al) = O(\a)
\]
holds, where the implicit constant depends only on~$\Om$.
\end{prp}

\begin{proof}
Let $!\p$ be the extension of Proposition~1. Then
\begin{align*}
m(\Om-\tilde\Om_\al)
&= m(\Om-\p(\Oma)) \\
&= m(\Om-!\p(\Oma)) \\
&= m(!\p(\Om-\Oma)) \\
&\le \abs{!\p:\lip,\Om} m(\Om-\Oma)  \\
&= O(\a).\qed
\end{align*}
\end{proof}

The Cantor family of torus embeddings may also be extended to a
family over~$\Om$. This extension can even be chosen so that the
additional embedded tori are still Lagrangian, though  of course
not invariant.

\begin{prp}
There exists an extension of $\Phi$ to a Lipschitz continuous
family of \ra. torus embeddings
\[
  !\Phi\maps \T^n\x\Om \to B\x\T^n
\]
such that each embedded torus is Lagrangian, and the estimates
for $!\Phi$ are the same as for~$\Phi$, though with a different
constant~$c$.
\end{prp}

\subdiv{d}

We thus arrive at the following conclusion.

\begin{proclaim}[Theorem B]
Suppose the assumptions of Theorem~A are satisfied. Then there
exist a lipeomorphism $!\p\maps\Om\to\Om$ close to the
identity and a family of torus embeddings
\[
  !\Phi\maps \T^n\x\Om \to \R^n\x\T^n
\]
close to $\Phi_0$ such that for every parameter value
\[
  \omt \in \tilde\Om_\al = !\p(\Oma)
\]
the \ham. \vf. $\rbar1{X_H}_{\omt}$ admits an invariant Lagrangian 
Kronecker torus $\Ts_\om = !\Phi(\T^n,\om)$, where $\om=!\p\inv(\omt)$.
Moreover, the estimates for $!\p$ and~$!\Phi$ are the same as for 
$\p$ and $\Phi$ in Theorem~A, though with a different constant~$c$, and
\[
  m(\Om-\tilde\Om_\al) = O(\a),
\]
where the implicit constant depends only on~$\Om$.
\end{proclaim}

We will see at the end of section~5 that the map $\p$ actually
can be assigned $\om$-derivatives of every order on the Cantor
set~$\Oma$. This may be formalized by introducing the
intrinsically defined notion of a differentiable function on an
arbitrary closed set~\cite{St,Wh}. The point is that -- due to the
Whitney extension theorem -- such functions can be extended to
functions on the whole space with the same differentiability
properties. The upshot is that there even exists an extension of
$\p$ to a $C^\infty$-function~$\bar\p$ on~$\Om$. The same
applies to~$\Phi$ and leads to the notion of smooth foliations of
invariant tori over Cantor sets~\cite{Po2}.

\subdiv{d}

We now prove the classical \kam. theorem.
Introducing the frequencies as parameters we wrote the \ham. as
$H=N+P$, where
\[
  P=P_h+P_{\e}
\]
is real analytic on $B\x\T^n\x\bar\Om$, $B$ some small ball
around the origin in~$\R^n$. Thus we can \emph{fix} some small
$h$ and~$s$, with $s^\nu\le h$, so that $P$ is real analytic on
the complex domain $D_{r,s}\x\W_h$ for all small~$r$, and so
that
\[
  \abs{P:r,s,h} 
  \le \abs{P_h:r,s,h} +\abs{P_{\smash{f_\ep}}:r,s,h}
  \le Mr^2+F\e,
\]
where $M$ is a bound on the hessian of~$h$ and
$F={\sup_{p,q,\e}}\abs{f_*(p,q,\e)}$.

To meet the smallness condition of Theorems A and~B, we
choose~$r$ by requiring that $Mr^2=F\e$ and arrive at the condition
\[
  2F\e \le \gm\a rs^\nu = \gm\a s^\nu\smash[b]{\sqrt{F\e\over M}},
\]
or
\[
  \e \le {\gm^2\a^2s^{2\nu}\over4FM} = \d\a^2, \qquad
  \d = {\gm^2s^{2\nu}\over4FM}.
\]
So there is a $\d$ depending  on $n$, $\tau$ and~$H$ such that
Theorems A and~B apply for $\e\le\d\a^2$.

By construction, an orbit $(I(t),\t(t))$ for the \ham.~$H$ at the
parameter value~$\omt$ translates into an orbit
$(p(t),q(t))=(p_0(\omt)+I(t),\t(t))$ for this \ham. in
$p,q$-\co.s.
It therefore follows with Theorem~B that the mapping
\[
  \Psi \maps \T^n\x\Om \to D\x\T^n,
\]
which is a composition of $!\Phi$ and $!\p$ with
\begin{align*}
  B\x\T^n\x\Om &\to D\x\T^n, \\
  (I,\t,\om) &\mapsto (h_p\inv(\om)+I,\t),
\end{align*}
is, for every $\om\in\Oma$, an embedding of an
invariant Lagrangian Kronecker torus $(\T^n,\om)$. Moreover, $\Psi$ is \emph{Lipschitz close} to the
real analytic unperturbed embedding
\begin{align*}
  \Psi_0\maps \T^n\x\Om &\to D\x\T^n, \\
  (\t,\om) &\mapsto (h_p\inv(\om),\t).
\end{align*}
It follows that the measure of the complement of all those tori
in the phase space is bounded by a constant times the measure of
$\Oma\x\T^n$, hence is~$\bigo\a$. This finishes the proof of the
classical \kam. theorem.

\section{Outline of the Proof of Theorem~A}

\subdiv{a}

We prove Theorem~A by a rapidly converging iteration procedure
that was proposed by Kolmogorov~\cite{Kol}. At each step
of this scheme a \ham. 
\[
  H_j = N_j+P_j
\] 
is considered, which is a small perturbation of  $N_j=e_j+\pair
\om I $. A \co. and parameter \trans. $\Fs_j$ is constructed such
that 
\[
  H_j\comp\Fs_j=N_{j+1}+P_{j+1}
\]
with another normal form $N_{j+1}$ and a much smaller error
term~$P_{j+1}$. Namely,
\[
  \n*{P_{j+1}} \le C \n*{P_j}^\k
\]
with some exponent $\k>1$. Repetition of this process leads to a
sequence of \trans.s $\Fs_0,\Fs_1,\dots$, whose products
\[
\Fsj = \Fs_0\comp\Fs_1\comp\dots\comp\Fs_{j-1}
\]
converge to an embedding of an invariant Kronecker torus.

In the meantime a number of other proofs have been given, for
example by formulating some generalized implicit function
theorem suited for small divisor problems~\cite{Ze76}, or by
referring to an implicit function theorem in tame Frechet
spaces~\cite{Bo}. Recently, Salamon and Zehnder~\cite{SZ} gave a proof that
avoids \co. \trans.s altogether and works in configuration
space. Also, Eliasson~\cite{El} described a way of using power series
expansions and majorant techniques in a very tricky way.

But here we stick to the traditional method of proof, as it
probably is the most transparent way to get to know the basic
techniques. They are indeed quite flexible and robust, and not
at all restricted to perturbations of integrable \ham. systems.
As we mentioned in the beginning, these techniques rather amount
to a \emph{strategy} of how to approach a large class of
perturbation problems.

\subdiv{b}

To describe one cycle of this iterative scheme in more detail we
now drop the subscript~$j$.

First, the perturbation~$P$ is approximated by some \ham.~$R$ by
truncating its Taylor series in~$I$ at first order and its Fourier
series in~$\t$ at some suitable high order~$K$. The
approximation error $P-R$ will be small, and we now consider the
\ham. $\H=N+R$ instead of $H=N+P$. The purpose of this
approximation will become clear later.

The \trans. $\Fs$ consists of a parameter dependent symplectic
change of \co.s~$\Phi$ and a change~$\p$ of the parameters
alone:
\[
  \Fs=(\Phi,\p) \maps 
  (I,\t,\om) \mapsto (\Phi(I,\t,\om),\p(\om)).
\]
Moreover, this \co. \trans. is of the form
\[
  \Phi \maps (I,\t,\om) \mapsto (U(I,\t,\om),V(\t,\om)),
\]
where $U$ is \emph{affine in~$I$}, and $V$ is independent of~$I$.
Such \trans.s~$\Fs$ form a group~$\Gs$ under composition.

We obtain $\Phi$ as the time-1-map of the flow $X_F^t$ of a \ham.
$F$, which is affine in~$I$. That is, $\Phi$ is the
time-1-shift determined by the equations of motion
\[
  \dot I = -F_\t(I,\t,\om), \quad \dot\t=F_I(\t,\om).
\]
Then $\Phi$ is symplectic for each~$\om$. 

To describe the transformed \ham. $H\comp\Phi$ we recall that
for a function~$K$,
\[
  \ddt\, K\comp X_F^t = \po KF \comp X_F^t,
\]
the \emph{Poisson bracket} of $K$ and~$F$ evaluated at~$X_F^t$.
Indeed,
\begin{align*}
  \rbar2{\ddt K\comp X_F^t}_{t=0} 
  &= \sum_{1\le j\le n} K_{\t_j}\dot\t_j + K_{I_j}\dot I_j \\
  &= \sum_{1\le j\le n} K_{\t_j}F_{I_j}-K_{I_j}F_{\t_j} 
   = \po KF,
\end{align*}
and the general formula follows.

So we can use Taylor's formula to expand $\H\comp\Phi =
\eval \H\comp X_F^t \at{t=1}$ \wrt. $t$ at~$0$ and write
\begin{align*}
  \H\comp\Phi
  &= \eval N\comp X_F^t\at{t=1} + \eval R\comp X_F^t \at{t=1} \\
  &= N + \po NF + \int_0^1 (1-t)\po{\po NF}F \comp X_F^t \dt \\
  &\qquad       + R + \int_0^1 \po RF \comp X_F^t \dt \\
  &= N+\po NF+R + \int_0^1 \po{(1-t)\po NF+R}F \comp X_F^t\dt.
\end{align*}
This is a \emph{linear} expression in $R$ and~$F$ -- the
\emph{linearization} of $\H\comp\Phi$ -- plus a \emph{quadratic}
integral remainder. That is, if $R$ and~$F$ are both roughly of
order~$\e$, then the integral will roughly be of order
$\e^2\ll\e$ and may be ascribed to the next perturbation
$P_\tp$.

The point is to find $F$ such that
\[
  N+\po NF + R = N_\tp
\]
is again a normal form. Equivalently, we want to solve
\[
  \po FN + \Nh = R, \qquad \Nh = N_\tp-N                   
  \eqlabel{LE}
\]
for $F$ and $\Nh$ when $R$ is given. Suppose for a moment that
such a solution exists. Then 
$(1-t)\po NF+R = (1-t)\Nh+tR$, and altogether we obtain
\[
  H\comp\Phi
  = \H\comp\Phi + (P-R)\comp\Phi = N_\tp + P_\tp
\]
with $N_\tp=N+\Nh$ and
\[
  P_\tp 
  = \int_0^1 \po{(1-t)\Nh+tR}F \comp X_F^t\dt + (P-R)\comp X_F^1                         
  \eqlabel{P+}
\]
as the new error term.

\subdiv{c}

Let us consider equation \eqref{LE} first on a formal level.
Clearly,
\[
  \del_\om F 
  \defeq \po FN = \sum_{1\le j\le n} F_{\t_j}N_{I_j}
  = \sum_{1\le j\le n} \om_j F_{\t_j}
\]
is a first order partial differential operator on the
torus~$\T^n$ with constant coefficients~$\om$. Expanding $F$ into
a Fourier series,
\[
  F = \sum_{k\in\Z^n} F_k e^{\i\pair k\t},
\]
with coefficients depending on $I$ and~$\om$, we find
\[
  \del_\om F = \sum_{k\in\Z^n} \i\!\pair k\om F_k e^{\i\pair k\t}.
\]
Thus, $\del_\om$ admits a basis of eigenfunctions $e^{\i\pair k\t}$
with eigenvalues $\i\!\pair k\om$, $k\in\Z^n$. In other words,
$\del_\om$ diagonalizes \wrt. this basis.

If $\om$ is now \emph{nonresonant}, then these eigenvalues are
all different from zero except when $k=0$. We then can solve for
all Fourier coefficients $R_k$ of the given function~$R$, except
for~$R_0$, which is given by the mean value of
$R$ over~$\T^n$,
\[
  R_0 = [R] = {1\over\operatorname{vol} \T^n} \int_{\T^n} R\dth.
\]
Hence, if $R$ is given, then we can always solve the equation $
\del_\om F = R-[R] $ at least formally by setting
\[
  F = \sum_{0\ne k\in\Z^n} {R_k\over \i\!\pair k\om}\,e^{\i\pair k\t}.                                                         
  \eqlabel{F}
\]
We are still free to add a $\t$-independent function to~$F$, but
we choose to normalize $F$ so that $[F]=0$.
-- 
Finally, equation \eqref{LE} is completely solved by setting
\[
  \Nh = [R].
\]
Of course, this choice of $\Nh$ is in no way uniquely determined,
but this is in some sense the simplest one.

\subdiv{d} 

There is a more systematic interpretation of the
preceding construction. For irrational~$\om$, the domain
of~$\del_\om$, consisting of all formal Fourier series in~$\t$
(ignoring the other \co.s here), \emph{splits} into two
\emph{invariant subspaces}, its nullspace~$\Ns$ consisting of
all constant functions, and its range~$\Rs$, consisting of all
series with vanishing constant term. Moreover, $\del_\om$ is
invertible on~$\Rs$.

Decompose $R$ into its respective components in $\Ns$ and~$\Rs$,
\[
  R = R_\Ns + R_\Rs.
\]
The projection onto $\Ns$ is given by taking the mean
value, so
\[
  R_\Ns = [R], \qquad R_\Rs = R-[R].
\]
The equation
\[
  \del_\om F + \Nh = R = R_\Rs + R_\Ns
\]
is then solved by ‘solving componentwise’,
\[
  \Nh = R_\Ns = [R], \qquad
  \del_\om F = R_\Rs = R-[R],
\]
where the latter can be solved uniquely for $F$ in~$\Rs$, since
$\del_\om$ is invertible on~$\Rs$.

This general procedure -- ‘solve for all the terms you can solve
for, and keep the rest’~-- 
is at the basis of all normal form theory. It just happens to take a
particularly simple form in our case.

\subdiv{e}

So far our considerations were formal. But in estimating the
series representation \eqref{F} of~$F$, we are confronted with the
well known and notorious problem of ‘small divisors’. Even if
$\om$ is nonresonant, infinitely many of the divisors
$\pair k\om $ become arbitrarily small in view of \eqref{Dir},
threatening to make the series \eqref{F} divergent.

This divergence is avoided, if $\om$ is required to be
strongly nonresonant. To formulate this key lemma, let
$\As^s$ denote the space of all analytic functions~$u$ defined
in the complex strip $\{\t:\sup_j\n*{\Im\t_j}<s\}
\subset \C^n$ with bounded sup-norm $\abs{u:s}$ over that strip.
Let \[
  \As^s_0 = \set{u\in\As^s: [u]=0},
\]
and recall that $\om\in\Dl_a^\tau$ satisfies
\[
  \abs{\pair k\om}\ge\a/\abs k^\tau \qt{for all} 0\ne k\in\Z^n.
\]

\begin{lem}
Suppose that $\om\in\Dl_\a^\tau$. Then the equation
\[
  \del_\om u = v, \qquad v\in\As^s_0,
\]
has a unique solution $u$ in $\bigcup_{0<\s<s} \As^{s-\s}_0$, with
\[
  \abs{u:s-\s} \le {c\over\a\s^{\tau}} \abs{v:s}, \qquad
  0<\s<s,
\]
where the constant $c$ depends only on $n$ and~$\tau$.
\end{lem}

\begin{proof}
We prove the lemma with $\s^{\tau+n}$ in place of $\s^\tau$ to
avoid lengthy technicalities. The interested reader is referred
to \cite{RuDS} or \cite{Po1} for a proof of the sharper result.
--
Expanding $u$ and $v$ into Fourier series, the unique formal
solution $u$ with $[u]=0$ is
\[
  u = \sum_{0\ne k\in\Z^n} {v_k\over \i\!\pair k\om}\,e^{\i\pair k\t}.
\]
As to the estimate we recall that the Fourier coefficients of an
analytic function on $\T^n$ decay exponentially fast:
\[
  \abs{v_k} \le \abs{v:s}e^{-\abs ks},
\]
where $\abs k = \abs{k_1}+\dots+\abs{k_n}$. See Lemma~A.1 for a
reminder. Together with the small divisor estimate for~$\om$ we
obtain
\[
  \abs{u:s-\s}
\le \sum_{k\ne0} {\abs{v_k}\over\abs{\pair k\om}}e^{\abs k(s-\s)}
\le {\abs{v:s}\over\a} \sum_{k\ne0} \abs k^\tau e^{-\abs k\s}.
\]
The infinite sum is now easily estimated by
constant times $\s^{-\tau-n}$.\qed 
\end{proof}

We observe that $\del_\om\inv$ is \emph{unbounded} as an
operator in $\As^s_0$. It is \emph{bounded} only as an operator
from $\As^s_0$ into the larger spaces $\As^{s-\s}_0$, with its
bound tending to infinity as $\s$ tends to zero. This phenomenon
is known as ‘loss of smoothness’ affected by the solution
operator $\del_\om\inv$, and is the main culprit why small
divisor problems are technically so involved. During the
iteration we have to let $\s\to0$ in order to stay in the
classes~$\As^s_0$. But then the
operator ${\del_\om\inv}$ is getting unbounded.  By the rapid
convergence of the Newton scheme, however, the error term
converges to zero even faster, thus allowing to overcome this
effect of the small divisors.

It is absolutely essential for Lemma~1 to be true that $\om$
satisfies infinitely many small divisor conditions, thus
restricting $\om$ to a \emph{Cantor} set with no interior points.
On the other hand, we will also need to transform the
frequencies and thus want them to live in \emph{open} domains.
This conflict is resolved by approximating $P$ by a trigonometric
polynomial~$R$. Then only finitely many Fourier coefficients need
to be considered at each step, and only \emph{finitely many} small
divisor conditions need to be required, which are easily
satisfied on some open $\om$-domain. Of course, during the
iteration more and more conditions have to be satisfied, and in
the end these domains will shrink to some Cantor set.

\subdiv{f}

We still have to finish one cycle of the iteration. Solving
\eqref{LE} we arrive at $H\comp\Phi=N_\tp+P_\tp$, where
\[
  N_\tp = N+[R] = e_\tp(\om)+\pair{\om+v(\om)}I,
\]
since $[R]$ is affine in~$I$ and independent of~$\t$. To
write $N_\tp$ again in normal form, we have to introduce
\[
  \om_\tp = \om+v(\om)
\]
as new frequencies. Since $v$ is small, there exists an inverse
map $\p\maps\om_\tp\mapsto\om$ by the implicit function
theorem -- see Appendix~A. With this change of parameters,
\[
  N_\tp = e_\tp(\om_\tp) + \pair{\om_\tp}I
\]
is again in normal form. This finishes one cycle of the
iteration.



\section{The \kam. Step}

\subdiv{a}

To avoid a flood of constants we will write
\[
  u\led v, \qquad u\dle v,
\]
if there exists a positive constant $c\ge1$, which depends
\emph{only on $n$ and~$\tau$} and could be made explicit, such
that $u \le cv$ and $ cu\le v$, respectively.
--
Now let $P$ be a real analytic perturbation of some normal
form~$N$.

\begin{proclaim}[The \kam. Step]
Suppose that $\abs{P:r,s,h}\le\e$ with
\[
  \tag{a}
  \e \dle \a\eta r\s^\nu, \qq
  \textup{(b)}\q  \e \dle hr, \qq
  \textup{(c)}\q  h \le {\a\over2K^\nu},
\]
for some $0<\eta<1/8$, $0<\s<s/5$ and $K\ge1$, where
$\nu=\tau+1$. Then there exists a real analytic \trans.
\[
  \Fs=(\Phi,\p) \maps 
  D_{\eta r,s-5\s}\x\W_{h/4} \to D_{r,s}\x\W_h
\]
in the group~$\Gs$ such that
$H\comp\Fs=N_\tp+P_\tp$ with
\[
  \abs{P_\tp:\eta r,s-5\s,h/4} \led {\e^2\over\a r\s^\nu} +
  \pas1{\eta^2+K^ne^{-K\s}}\e.
\]
Moreover, on $D_{\eta r,s-5\s}\x\W_{h}$ and $\W_{h/4}$ the estimates
\begin{align*}
  \abs{W(\Phi-\id)}, \n*{W(D\Phi-\Id)W\inv} 
  &\led {\e\over\a r\s^\nu}, 
  \\
  \abs{\p-\id}, h\abs{D\p-\Id} &\led {\e\over r},
\end{align*}
hold with the weight matrix $W=\diag(r\inv\Id,\s\inv\Id)$.
\end{proclaim}

\subdiv{b}

The proof of the \kam. Step follows the lines of the
preceding section and consists of six small steps. Except for
the last step everything is uniform in $\W_h$, whence we write
$\abs{\cd:r,s}$ for $\abs{\cd:r,s,h}$ throughout.

\paragraph{1. Truncation}
We approximate $P$ by a \ham.~$R$, which is affine in~$I$ and a
trigonometric polynomial in~$\t$. To this end, let $Q$ be the
linearization of $P$ in~$I$ at $I=0$. By Taylor's formula with
remainder and Cauchy's estimate -- see Appendix~A for a reminder~--,
we have 
\[
  \abs{Q:r,s} \led \e, \qquad \abs{P-Q:2\eta r,s} \led \eta^2\e.
\]
Then we simply truncate the Fourier series of~$Q$ at order~$K$
to obtain~$R$. By Lemma~A.2,
\[
  \abs{R-Q:r,s-\s} \led K^ne^{-K\s}\e.
\]
Since the factor $K^ne^{-K\s}$ will be made small later on, we
also have
\[
  \abs{R:r,s-\s} \led \e.
\]
See Appendix~A for some remarks about this truncation of Fourier
series.

\paragraph{2. Extending the small divisor estimate}
The nonresonance conditions \eqref{Dio} are assumed to hold on
the real set $\Oma$ only. But assumption~(c) implies that 
\[
  \abs{\pair k\om} \ge {\a\over2\abs k^\tau} \quad
  \qtext{for all} 0\ne\abs{k}\le K
  \eqlabel{NR3}
\]
for all $\om$ in the \nbhd. $\W_h$ of~$\Oma$. Indeed, for every
$\om\in\W_h$ there is some $\om_*\in\Oma$ with
$\abs{\om-\om_*}<h$, hence
\[
  \abs{\pair k{\om-\om_*}} 
  \le \abs k \abs{\om-\om_*}
  \le Kh \le {\a\over2K^\tau} 
  \le {\a\over2\abs k^\tau}
\]
for $\abs k\le K$. Together with the estimate \eqref{Dio} for
$\pair k{\om_*}$ this proves the claim.

\paragraph{3. Solving the linearized equation $\protect\po FN+\Nh=R$}
We solve this equation as described in the preceding section. We
have $\Nh=[R]$ and thus
\[
  \abs{\Nh:r} \le \abs{R:r,s-\s} \led \e.
\]
We can solve for $F$ uniformly for all~$\om$ in $\W_h$ because of
\eqref{NR3} and the fact that $R$ contains only Fourier
coefficients up to order~$K$, by truncation. Hence the estimate
of Lemma~1 applies as well, and we obtain a real
analytic function~$F$ with
\[
  \abs{F:r,s-2\s} 
  \led {\abs{R:r,s-\s}\over \a\s^{\tau}}
  \led {\e\over\a\s^{\tau}}.
\]
With Cauchy we get
$\abs{F_\t:r,s-3\s} \led \e/\a\s^{\tau+1}$ and
$\abs{F_I:r/2,s-2\s} \led \e/\a r\s^{\tau}$, hence
\[
  \frac1r\abs{F_\t},
  \frac1\s\abs{F_I} \led {\e\over\a r\s^\nu}
\]
uniformly on $D_{r/2,s-3\s}$ with $\nu=\tau+1$.

\paragraph{4. Transforming the \co.s}
The \co. \trans. $\Phi$ is obtained as the real
analytic time-1-map of the flow $X_F^t$ of the \ham.
\vf.~$X_F$, with equations of motions
\[
  \dot I = -F_\t, \qquad \dot\t=F_I.
\]
With assumption~(a) and the preceding estimates we can assure
that we have $\abs{F_\t}\le\eta r\le r/8$ and $\abs{F_I}\le\s$ on
$D_{r/2,s-3\s}$ uniformly in~$\om$. Therefore, the time-1-map
is well defined on $D_{r/4,s-4\s}$, with
\[
  \Phi=\eval X_F^t\at{t=1} \maps D_{r/4,s-4\s} \to D_{r/2,s-3\s},                                                
  \eqlabel{Phimap}
\]
and
\begin{align*}
  \abs{U-\id} &\le \abs{F_\t} \led {\e\over\a\s^\nu} , 
  \\
  \abs{V-\id} &\le \abs{F_I} \led {\e\over\a r\s^{\nu-1}}
\end{align*}
on that domain for $\Phi=(U,V)$.
The Jacobian of $\Phi$ is
\[
  D\Phi = \pmat{ U_I&U_\t\\ 0 &V_\t} ,
\]
since $F$ is \emph{linear} in~$I$, hence $F_I$ and $V$ are
\emph{independent} of~$I$. By the preceding estimates and
Cauchy,
\[
  \abs{U_I-\Id} \led {\e\over\a r\s^\nu}, \qquad
  \abs{U_\t} \led {\e\over\a\s^{\nu+1}}, \qquad
  \abs{V_\t-\Id} \led {\e\over\a r\s^\nu}
\]
on the domain $D_{r/8,s-5\s} \supseteq D_{\eta r,s-5\s}$. This
proves the estimates for~$\Phi$. Finally, we observe that
$\abs{U-\id} \le\abs{F_\t} \le \eta r$ implies that also
$\Phi \maps D_{\eta r,s-5\s} \to D_{2\eta r,s-4\s}$.

\paragraph{5. New error term}
To estimate $P_\tp$ as given in \eqref{P+} we first consider the term
$\po RF$. Again, by Cauchy's estimate,
\begin{align*}
  \abs{\po RF:r/2,s-3\s}
  &\led \abs{R_I}\abs{F_\t}+\abs{R_\t}\abs{F_I}  \vphantom{\frac.r}\\
  &\led \frac\e r\cdot \frac\e{\a\s^\nu} +
             \frac\e\s\cdot\frac\e{\a r\s^{\nu-1}} \\
  &\led {\e^2\over\a r\s^\nu}.
\end{align*}
The same holds for $\abs{\po\Nh F:r/2,s-3\s}$. Together with
\eqref{Phimap} and $\eta<\frac18$ we get
\begin{multline*}
  \abs{\vphantom\sum\smash{\int_0^1 \po{(1-t)\Nh+tR}F \comp X_F^t\dt} :\eta r,s-5\s} \\
  \le \abs{\po{(1-t)\Nh+tR}F:r/2,s-4\s} 
  \led {\e^2\over\a r\s^\nu}.
\end{multline*}
The other term in \eqref{P+} is bounded by 
\begin{align*}
  \abs{(P-R)\comp\Phi:\eta r,s-5\s}
  &\le \abs{P-R:2\eta r,s-4\s} \\
  &\le \abs{P-Q:2\eta r,s-4\s} +\abs{Q-R:2\eta r,s-4\s} \\
  &\led \pas1{\eta^2+K^ne^{-K\s}}\e.
\end{align*}
These two estimates together give the bound for~$\abs{P_\tp}$.

\paragraph{6. Transforming the frequencies}
Finally, we have to invert the map
\[
  \om \mapsto \om_\tp=\om+v(\om), \qquad  v=\Nh_{\!I}=[R_I]
\]
to put $N+\Nh$ back into a normal form~$N_\tp$. With
assumption~(b) and Cauchy's estimate we can assure that
\[
  \abs{v:h/2} = \abs{\Nh_{\!I}:h/2} \led \frac\e r \le \frac h4.
\]
The implicit function theorem of Appendix~A applies, and there
exists a real analytic inverse map
$ \p \maps \W_{h/4} \to \W_{h/2}$, $ \om_\tp\mapsto\om $, with
the estimates
\[
  \abs{\p-\id}, h\abs{D\p-\Id} \led \frac\e r
\]
on $\W_{h/4}$.  Setting
$N_\tp=(N+\Nh)\comp\p$ the proof of the \kam. Step is completed.

\section{Iteration and Proof of Theorem A}

\subdiv{a}

\def\cc{c}%
We are going to iterate the \kam. Step infinitely often,
choosing appropriate sequences for the parameters $\s$, $\eta$
and so on. To motivate our choices, let us start by fixing a
geometric sequence for~$\s$, say, $\s_\tp=\s/2$, where the plus
sign indicates the corresponding parameter value for the next
step. Let $r_\tp=\eta r$, and let us consider the \emph{weighted
error terms}
\[
  E={\e\over\a r\s^\nu}, \qquad E_\tp={\e_\tp\over\a r_\tp\s_\tp^\nu}.
\]
Then we have
\[
  E_\tp \led {E^2\over\eta}+\pas1{\eta^2+K^ne^{-K\s}}{E\over\eta}.
\]
Suppose we can choose $\eta$ and~$K$ so that $\eta^2=E$ and
$K^ne^{-K\s}\le E$. Then
\[
  E_\tp \led \eta\inv E^2 = E^\k, \qquad \k=\frac32.
\]
That is, $E_\tp\le \cc^{\k-1}E^\k$ for some constant~$\cc$
determined by the \kam. Step and depending only on $n$
and~$\tau$. Consequently,
\[
  \cc E_\tp \le \cc^\k E^\k,
\]
and this scheme converges exponentially fast for $ E<\cc\inv$.

It remains to discuss our assumptions
\begin{align*}
  \qquad \eta^2 &= E  \tag d 
  \\
  K^ne^{-K\s} &\le E  \tag e
\end{align*}
as well as assumptions (a--c) of the \kam. Step. There is no
obstacle to take (d) as the \emph{definition} of~$\eta$, as this
\emph{implies} $E\dle\eta$ and thus (a) for $E$ \suff. small. The
other three conditions amount to
\[
  {\e\over r} \dle h \le {\a\over2K^\nu}, 
  \qquad
  K^ne^{-K\s} \le E = {\e\over\a r\s^\nu}.
\]
They are combined into the sufficient condition
\[
  K^\nu\s^\nu e^{-K\s} 
  \le E\s^\nu 
  = {\e\over\a r}
  \dle {h\over\a} 
  \le {1\over2K^\nu}. 
  \eqlabel{f}
\]
Now it suffices to set up geometric sequences for $K$ and~$h$ such
that, say,
\[
  K_\tp \s_\tp = 2K\s, \qquad K_\tp^\nu h_\tp \le K^\nu h.
\]
Then these inequalities hold inductively, provided they hold
initially, with $K_0\s_0$ \suff. large. In particular, after
fixing $\s_0$, $K_0$ and~$E_0$, we may set
\[
  h_0 =\a c_0E_0\s_0^\nu
\]
with a suitable constant~$c_0$.

\subdiv{b}

We are now ready to set up our parameter sequences. Let
\[
  \s_{j+1}={\s_j\over2}, \qquad
  K_{j+1} = 4K_j, \qquad
  h_{j+1} = {h_j\over4^\nu},
\]
where $\s_0=s_0/20$, and $K_0$ is chosen so large that the left
hand side of~(f) is smaller than its right hand side, and small
enough so that $E_0$ can be fixed to meet (a) and~(d). In
addition, $K_0\s_0$ has to be so large that the left hand side
of~(f) decreases at least at an exponential rate~$\k$. Thus we
need $1\dle K_0\s_0$. Subsequently, we fix $h_0$ as above. 

Next, let
\[
  E_{j+1} = \cc^{\k-1}E_j^\k, \qquad
r_{j+1} = \eta_j r_j, \qquad
\eta_j^2 = {E_j},
\]
where $r_0$ is still free and $\cc$ given by the \kam. Step.
Finally, we define 
\[
  s_{j+1}=s_j-5\s_j
\]
and the complex domains
\begin{align*}
  D_j  &= \{\abs I<r_j\} \x \{\abs{\Im\t}<s_j\}, \\
  \W_j &= \{\abs{\om-\Oma}<h_j\}.
\end{align*}
Note that $s_j \to s/2$ and $r_j\to 0$. --
Now let $H=N+P_0$.

\begin{proclaim}[Iterative Lemma]
Suppose $P_0$ is real analytic on $D_0\x\W_0$ with
\[
  \abs{P_0:r_0,s_0,h_0} \le \e_0 \defeq \a E_0 r_0\s_0^\nu.
\]
Then for each $j\ge0$ there exists a normal form~$N_j$ and a real
analytic \trans. 
\[
  \Fsj=\Fs_0\comp\dots\comp\Fs_{j-1}\maps D_j\x\W_j\to D_0\x\W_0
\]
in the group~$\Gs$ such that $H\comp\Fsj=N_j+P_j$ with
\[
  \abs{\smash{P_j}:r_j,s_j,h_j} \le \e_j \defeq \a E_jr_j\s_j^\nu.
\]
Moreover,
\[
  \abs{!W_0(\Fs^{\j+1}-\Fsj)} \led {\e_j\over r_jh_j}
\]
on $D_{j+1}\x\W_{j+1}$ with the weight matrix
$!W_0 = \diag(r_0\inv\Id,\s_0\inv\Id,h_0\inv\Id)$.
\end{proclaim}

\begin{proof} 
Letting $\Fs_0=\id$, there is nothing to do
for $j=0$. To proceed by induction, we have to check the
assumptions of the \kam. Step for each $j\ge0$. But (a) is
satisfied by the definition of $\eta_j$ and the \suff. small
choice of~$E_0$, and (b--c) hold by the definition of $h_j$ and
$K_j$ and the choice of their initial values.

We obtain a transformation
\[
  \Fs_j\maps D_{j+1}\x\W_{j+1}\to D_j\x\W_j
\] 
taking $H_j=N_j+P_j$
into $H_j\comp\Fs_j=N_{j+1}+P_{j+1}$ with
\begin{multline*}
  \n*{P_{j+1}} 
   \led \e_jE_j+\pas1{\eta_j^2+K_j^ne^{-K_j\s_j}}\e_j \\
   \led \e_jE_j 
   = \a E_j^2r_j\s_j^\nu 
   \led \a \eta_j\inv E_j^2r_{j+1}\s_{j+1}^\nu.
\end{multline*}
Since $\eta_j\inv E_j^2 = E_j^\k = \cc^{1-\k}E_{j+1}$,  we obtain 
$\n*{P_{j+1}}\le\e_{j+1}$ by an appropriate choice of~$\cc$ as
required. Thus, the transformation
$\Fs^{\j+1}=\Fsj\comp\Fs_j=\Fs_0\comp\dots\comp\Fs_j$  takes $H$ into
$N_{j+1}+P_{j+1}$ with the proper estimate for~$P_{j+1}$.

The estimate of $\Fsj$ requires a bit more, though elementary
work. We observe that the estimates of the \kam. Step and
Cauchy imply
\[
  \abs{!W_j(\Fs_j-\id)},
  \abs{!W_j(!D\Fs_j-\Id)!W_j\inv}
  \led \max\pas3{{\e_j\over\a r_j\s_j^\nu},{\e_j\over r_jh_j}}
  \led {\e_j\over r_jh_j}
\]
on $D_{j+1}\x\W_{j+1}$, where
$!D$ denotes the Jacobian \wrt. $I$, $\t$ \emph{and}~$\om$, and
$ !W_j = \diag(r_j\inv\Id,\s_j\inv\Id,h_j\inv\Id) $. We then have
\begin{multline*}
  \abs{!W_0(\Fs^{\j+1}-\Fsj)}
   = \abs{ !W_0(\Fsj\comp\Fs_j-\Fsj) } \\
   \le \abs{ !W_0!D\Fsj!W_j\inv} \abs{!W_j(\Fs_j-\id)} 
   \led \abs{ !W_j(\Fs_j-\id)} 
   \led {\e_j\over r_jh_j} \kern-1pc
\end{multline*}
\emph{provided} we can bound the first factor in the
second row on the domain $D_j\x\W_j$.

But by induction we have 
$!D\Fsj=!D\Fs_0\comp\dots\comp!D\Fs_{j-1}$, with the Jacobians
evaluated at different points, which we do not indicate. Since 
$\n*{!W_i!W_{i+1}\inv}\le1$ for all~$i$, we can use a telescoping
argument and the inductive estimates for the $\Fs_i$ to obtain 
\begin{align*} 
  \abs{!W_0!D\Fsj!W_j\inv}
  &\le \abs{ !W_0!D\Fs_0\comp\dots\comp\Fs_{j-1}!W_j\inv} \\[3pt]
  &\le \abs{!W_0!D\Fs_0!W_0\inv} \abs{!W_0!W_1\inv} \x  \\
  & \qquad \cdots \x \abs{!W_{j-1}!D\Fs_{j-1}!W_{j-1}\inv} \abs{!W_{j-1}!W_{j}\inv} \\
  &\le \prod_j \pas3{1+{c_1\e_j\over r_jh_j}}.
\end{align*}
This is uniformly bounded and small,
since $\e_j/r_jh_j$ converges rapidly to zero.
\end{proof}

\subdiv{c}

We prove Theorem~A by applying the Iterative Lemma to the hamiltonian
$H=N+P$, letting $P_0=P$, $r_0=r$, $s_0=s$. We have
$h_0 \dle \a\s_0^\nu \le \a s^\nu \le h$ by construction and
assumption, and 
\[
  \abs{P_0:r_0,s_0,h_0} \le \abs{P:r,s,h} \le \e
\le \gm\a rs \le \e_0 = \a E_0r_0\s_0^\nu,
\]
by fixing the constant $\gm$ in Theorem~A \suff. small.

By the estimates of the Iterative Lemma the maps $\Fsj$ converge
uniformly on
\[
  \bigcap_{j\ge0} D_j\x\W_j = T_*\x\Oma, \qquad
T_* = \{0\}\x\{\abs{\Im\t}<s/2\},
\]
to a map~$\Fs$ consisting of a family of embeddings
$\Phi\maps \T^n\x\Oma \to D\x\T^n$ and a parameter
transformation $\p\maps\Oma\to\Om$, which are \ra. on~$\T^n$
and uniformly continuous on~$\Oma$. Moreover,
\[
  \abs{!W_0(\Fs-\id)} \led {\e_0\over r_0h_0}
\]
uniformly on $T_*\x\Oma$ by the usual telescoping argument.

From the estimate $\n*{H\comp\Fsj-N_j}\led\e_j$ on $D_j\x\W_j$
we obtain 
\[
  \abs{W_j(J(D\Phi^j)^t\nabla H\comp\Fsj - J\nabla N_j)}
  \led {\e_j\over r_j\s_j}
\]
on $T_*\x\W_j$ with $W_j=\diag(r_j\inv\Id,\s_j\inv\Id)$.
The symplectic nature of the map $\Phi^j$ and the uniform
estimate of $!W_0!D\Fsj!W_j\inv$ above then imply
\[
  \abs{X_H\comp\Fsj - D\Phi^j{\cdot}X_N} \led {\e_j\over r_j\s_j}
\]
on $T_*\x\Oma$ for all~$j$, where $X_N$ is the \ham. \vf.
of $N=\pair\om I$. Going to the limit we obtain
\[
  X_H\comp\Fs = D\Phi{\cdot}X_N.
\]
Thus, $\Phi$ is an embedding of the Kronecker torus $(\T^n,\om)$
as an invariant torus of the \ham. \vf.~$X_H$ at the
parameter~$\p(\om)$. Moreover, this torus is Lagrangian, since
\[
  \Phi^*\ups  
  = \lim_{j\to\infty} \eval(\Phi^j)^*\ups\at{I=0} 
  = \lim_{j\to\infty} \eval\ups\at{I=0} = 0
\]
by the symplecticity of the~$\Phi^j$.

\subdiv{d}

Let us now look at the $\om$-derivatives of the $\Fsj$. Since
$E_j$ converges to zero at an \emph{exponential} rate, we have
\[
  {\e_j\over r_jh_j^m} \to 0 \q\qtext{for all} m\ge0.
\]
Hence, \emph{all} $\om$-derivatives of the $\Fsj$ converge
uniformly on $T_*\x\Oma$, and we could assign
$\om$-derivatives of any order to the limit map~$\Fs$ on the Cantor
set~$\Oma$~\cite{St}. Without making this concept precise,
however, we can at least conclude that $\Fs$ is \emph{Lipschitz
continuous} in~$\om$. Its Lipschitz norm is bounded by the limit
of the bounds on the first $\om$-derivatives of the
$\Fsj$. On $T_*\x\Oma$, the usual Cauchy estimate yields
\[
  \abs{!W_0(\Fs-\id):\lip} \led {\e_0\over r_0h_0^2}.
\]

\subdiv{e}

We finally look at the estimate of~$\Fs$. So far it does not
reflect the actual size~$\e$ of the perturbation, since we fixed
$E_0$ and thus $\e_0$ independently of~$\e$. But we observe that
everything is still all right if in all the estimates for $P_j$,
$\Fs_j$ and~$\Fsj$, the $\e_j$ are scaled down by the linear factor
\[
  {\e\over\e_0} =\a E_0r_0\s_0^\nu.
\]
Scaling down our estimates of~$\Fs$ by this factor we can finally
extract our estimates of $\Phi$ and~$\p$ as stated in
Theorem~A, since $\a\s_0^\nu\dle h_0$. This finishes the proof.
\vglue1.5pc
\noindent

\appendix

\section{Some Facts about Analytic Functions}

\subdiv{a}

First we recall a variant of the Cauchy estimate, which is 
used over and over. Let $D$ be an open domain in $\C^n$, let
$D_r=\set{z: \abs{z-D}<r}$ be the \nbhd. of radius~$r$
around~$D$, and let $F$ be an analytic function on~$D_r$ with
bounded sup-norm $\abs{f:r}$. Then
\[
  \n1{f_{z_j}}_{r-\r} \le \frac1\r\abs{f:r}
\]
for all $0<\r<r$ and $1\le j\le n$. This follows immediately from
the Cauchy estimate for one complex variable.

\subdiv{b}

Next we give the estimate for the Fourier coefficients of an
analytic function~$v$ on~$\T^n$ used in the proof of Lemma~1.
Recall that $\As^s$ denotes the space of all functions
on $\T^n$ bounded and analytic in the strip
$\{\t:\abs{\Im\t}<s\}$.

\begin{alm}
If $v\in\As^s$, then $v=\sum_k v_ke^{\i\pair k\t}$ with
\[
  \abs{v_k} \le \abs{v:s}e^{-\abs k s}, \qquad
  k\in\Z^n.
\]
\end{alm}

\begin{proof}
The Fourier coefficients $v_k$ of~$v$ are given by
\[
  v_k = {1\over(2\pi)^n} \int_{\T^n} v(\t)e^{-\i\pair k\t}\dth.
\]
Since the integral of an analytic function over a closed
contractible loop in any of the coordinate planes is zero, and
since $v$ is $2\pi$-periodic in each argument also in the
complex \nbhd., the path of integration may be shifted into the
complex, so that
\[
  v_k = {1\over(2\pi)^n} \int_{\T^n}
     v(\t-i\p)e^{-\i\pair k{\t-\i\p}}\dth
\]
for any constant real vector~$\p$ with $\abs{\p}<s$.
Choosing $\p=(s-\s)(e_1,\dots,e_n)$ with
$0<\s<s$ and $e_j = \text{sgn}\,k_j$, $1\le j\le n$, we obtain
\[
  \abs{v_k} \le \abs{v:s}e^{-\abs k(s-\s)}
\]
for all $\s>0$. Letting $\s\to0$ the lemma follows. \qed
\end{proof}

We can now also estimate very roughly the remainder, when we
truncate the Fourier series of~$v$ at order~$K$ to obtain
$T_Kv = \sum_{\abs k\le K} v_k e^{\i\pair k\t}$.

\begin{alm}
If $v\in\As^s$ and 
\footnote{The condition $K\sg\ge1$ is missing in the original version 
of this paper that appeared in print. I am thankful to San Vu Ngoc for pointing this out to me.}
$K\s\ge1$, then
\[
  \abs{v-T_Kv:s-\s} \le CK^ne^{-K\s}\abs{v:s}, \qquad 0\le\s\le s,
\]
where the constant $C$ only depends on~$n$.
\end{alm}

\begin{proof}
With Lemma~A.1,
\begin{align*}
  \abs{v-T_Kv:s-\s} 
  &\le \sum_{\abs k>K} \abs{v_k}e^{\abs k(s-\s)} \\
  &\le \abs{v:s} \sum_{\abs k>K} e^{-\abs k\s}
   \le \abs{v:s} \sum_{l>K} 4^n l^{n-1}e^{-l\s},
\end{align*}
by summing first over all $k$ with $\abs k=l$, whose number is
bounded by $4^nl^{n-1}$. The last sum is then
easily bounded by a constant times $K^ne^{-K\s}$. \qed
\end{proof}

There are much more efficient ways to approximate a periodic
function~$v$ by trigonometric polynomials. The above crude way
amounts to multiplying the Fourier transform $\hat v$ of~$v$
with a \emph{discontinuous} cut off function. Instead, one
should multiply $\hat v$ with a \emph{smooth} cut off
function~$\psi_K$. For instance, one could take $\psi_K(x) =
\psi(x/K)$, where $\psi$ is a fixed function, which is $1$ on the
ball $\abs x\le\frac12$, vanishes outside the ball $\abs x\ge1$, and
between $0$ and $1$ otherwise. Transforming back,
\[
  (\hat v \psi_K)\hat{\phantom{x}} = v \ast \hat\psi_K
\]
amounts to a convolution of $v$ with a real analytic
approximation of the identity $\hat\psi_K$, as $K\to\infty$.
Such smoothing operators have many interesting properties. For
more details, see for example~\cite{Ze76}.

\subdiv{c}

We finally formulate a version of the implicit function
theorem for analytic maps, which we need to invert the frequency
map during the \kam. Step. 
Recall that $\W_h$ is an open complex \nbhd. of radius~$h$
of some subset $\Om$ of $\R^n$.
In the following, $\abs\cd$ denotes the sup-norm for vectors
and maps, and the induced operator-norm for Jacobians.
 
\begin{alm}
Suppose $f$ is a real analytic map from $\W_h$ into~$\C^n$. If
\[
  \abs{f-\id:\W_h} \leq \d \leq h/4,
\]
then $f$ has a real analytic inverse $\p$ on
$\W_{h/4}$. Moreover,
\[
    \abs{\p-\id},\quad
  \frac h4\,\abs{D\p-\Id} \leq \d
\]
on this domain.
\end{alm}

\begin{proof}
Let $\eta=h/4$. Let $u,v$ be two points in $\W_{2\eta}$ with
$f(u)=f(v)$. Then
\[
  u-v = \pas{u-f(u)}-\pas{v-f(v)},
\]
hence
$ \abs{u-v} \leq 2\d \leq 2\eta $.
It follows that the segment $(1-s)u+sv$, $0\leq s\leq1$,
is strictly contained in $\W_{3\eta}$. Along this segment,
\[
  \t = \max\abs{D f-I} < \d/\eta \leq 1
\]
by Cauchy's inequality and so
\[
  \abs{u-v} 
  \leq \abs{D f-I}\abs{u-v}
  \leq \t\abs{u-v}
\]
by the mean value theorem. It follows that $u=v$.
Thus, $f$~is one-to-one on the domain~$\W_{2\eta}$.
 
By elementary arguments from degree theory the image of
$\W_{2\eta}$ under~$f$ covers $\W_\eta$, since
$\abs{f-\id}\leq\d$. So $f$ has a real analytic inverse $\p$ on
$\W_\eta$, which clearly satisfies $ \abs{\p-\id} \leq \d $.
Finally, 
\begin{align*}
  \abs{D\p-I}_{\eta}
  &\leq \n*{\pas{D f}\inv - I}_{2\eta}  \\
  &\leq {1\over1-\n*{D f - I}_{2\eta}} -1   
   \leq {1\over1-\d/2\eta}-1 
   \leq {\d\over \eta}
\end{align*}
by applying Cauchy to the domain~$\W_{2\eta}$.
\qed\end{proof}

\section{Lipschitz Functions}

Let $B\subset\R^n$ be a closed set. We prove the basic fact --
used in section~2 -- that a Lipschitz continuous function
$u\maps B\to\R$ can be extended to a Lipschitz continuous
function $U\maps\R^n\to\R$ without affecting its Lipschitz
constant
\[
  \abs{u:\lip,B} = \sup_{\substack{x,y\in B\\x\ne y}}
  {\abs{u(x)-u(y)}\over\abs{x-y}},
\]
where on $\R^n$ we may take any norm~$\abs\cd$. That is, we
have
\[
  \eval U\at B = u, \qquad \abs{U:\lip,\R^n} = \abs{u:\lip,B},
\]
Incidentally, for this extension $B$ could be \emph{any} point
set.

Indeed, $U$ is given by
\[
  U(x) = \sup_{z\in B} \pas{u(z)-\l\abs{z-x}}, \qquad x\in\R^n,
\]
where $\l=\abs{u:\lip,B}$. By the triangle inequality,
\[
  \pas{u(z)-\l\abs{z-y}} \ge \pas{u(z)-\l\abs{z-x}} - \l\abs{x-y}.
\]
Taking suprema over~$z$ on both sides we obtain $U(y)\ge U(x)-\l\abs{x-y}$, or
equivalently $ U(x)-U(y) \le \l\abs{x-y} $.
Interchanging $x$ and~$y$ we obtain
\[
  \abs{U(x)-U(y)} \le \l\abs{x-y},
\]
whence $\abs{U:\lip,\R^n}\le \abs{u:\lip,B}$. We leave it to
the reader to check that $U=u$ on~$B$.


\begin{thebibliography}{\hskip12.5pt}

\addcontentsline{toc}{section}{References}
\small
\itemsep1.5pt plus .2pt
\raggedright
\catcode`\…\active \def…{.\thinspace\ignorespaces}


\bibitem{A63}
\sc V…I… Arnol'd,
\rm Proof of a theorem by A.N. Kolmogorov on the invariance
of quasi-periodic motions under small perturbations of the
Hamiltonian.  
\it Russian Math. Surveys \bf{18} \rm(1963) 9--36.

\bibitem{AMM}
\sc V…I… Arnol'd,
\it Mathematical Methods of Classical Mechanics.
\rm Springer, 1978.

\bibitem{AEMS}
\sc V…I… Arnol'd {\rm (ed)},
\it Dynamical Systems III. 
\rm Encyclopaedia of Mathematical Sciences Volume 3, Springer, 1988.

\bibitem{Bo}
\sc J…-B… Bost,
\rm Tore invariants des systèmes dynamiques Hamiltoniens
(d'apres Kolmogorov, Arnold, Moser, Rüssmann, Zehnder, Herman,
Pöschel, \dots). 
\it Asté\-ris\-que \bf{133--134} \rm(1986) 113--157.

\bibitem{El}
\sc L…H… Eliasson,
\rm Absolutely convergent series expansions for quasi periodic 
motions. 
\it Math. Phys. Electron. J. \bf{2} \rm(1996) paper 4, 33 p. 

\bibitem{Ge}
\sc C… Genecand,
\rm Transversal homoclinic points near elliptic fixed
points of area-preserv\-ing diffeomorphisms of the plane.
\rm In \it Dynamics Reported  New Series, \rm Volume~2, Springer, 1993, 1--30.

\bibitem{Kat}
\sc A…B… Katok,
\rm Ergodic properties of degenerate integrable Hamiltonian
systems. 
\it Math. USSR Izv. \bf{7} \rm(1973) 185--214.

\bibitem{Kol}
\sc A…N… Kolmogorov,
\rm On the conservation of conditionally periodic motions
for a small change in Hamilton's function 
\it Dokl. Akad. Nauk SSSR \bf{98} \rm(1954) 527--530  [Russian]. 
English translation in \it Lectures Notes in Physics \bf 93\rm, Springer, 1979.

\bibitem{MM}
\sc L… Markus \& K… Meyer,
\rm Generic Hamiltonian systems are neither integrable nor
ergodic.
\it Mem. Am. Math. Soc. \bf{144} \rm(1974).

\bibitem{Mo63}
\sc J…Moser,
\rm On invariant curves of area preserving mappings of an
annulus.
\it Nachr. Akad. Wiss. Gött.,  Math. Phys. Kl. \rm (1962) 1--20.

\bibitem{Mo67}
\sc J…Moser,
\rm Convergent series expansions for quasi-periodic motions.
\it Math. Ann. \bf{169} \rm(1967) 136--176.

\bibitem{Mo69}
\sc J…Moser,
\rm On the continuation of almost periodic solutions for
ordinary differential equations.
\rm In \it Proc. Int. Conf. Func.
Anal. Rel. Topics, \rm Tokyo, 1969, 60--67.

\bibitem{MoR}
\sc J…Moser,
\it Stable and Random Motions in Dynamical Systems.
\rm Princeton University Press, 1973.

\bibitem{MZ}
\sc J…Moser \& E…Zehnder,
\it Notes on Dynamical Systems. 
\rm Courant Lecture Notes in Mathematics 12, New York, 2005.

\bibitem{Po1}
\sc J…Pöschel,
\rm Über invariante Tori in differenzierbaren Hamiltonschen
Systemen. 
\it Bonn. Math. Schr. \bf{120} \rm(1980) 1--103.

\bibitem{Po2}
\sc J…Pöschel,
\rm Integrability of Hamiltonian systems on Cantor sets.
\it Comm. Pure Appl. Math. \bf{35} \rm(1982) 653--695.

\bibitem{RuDS}
\sc H…Rüssmann,
\rm On optimal estimates for the solution of linear partial
differential equations of first order with constant coefficients
on the torus.
\rm In \sc J…Moser \rm (ed), \it Dynamical Systems, Theory and Applications,
\rm  Lecture Notes in Physics 38, Springer, 1975, 598--624.

\bibitem{RuDV}
\sc H…Rüssmann,
\rm Non-degeneracy in the perturbation theory of integrable
dynamical systems.  
\rm In \sc  M.M. Dodson \& J.A.G. Vickers \rm(eds),
\it Number Theory and Dynamical Systems,
\rm London Math. Soc. Lecture Note Series 134, Cambridge University
Press, 1989, 5--18.

\bibitem{Ru98}
\sc H…Rüssmann,
\rm Invariant tori in the perturbation theory of weakly
nondegenerate integrable Hamiltonian systems.
\rm Preprint, Mainz (1998).

\bibitem{Ru01}
\sc H…Rüssmann,
\rm Invariant tori in non-degenerate nearly integrable Hamiltonian systems.
\it Regular Chaotic Dynam. \bf 6 \rm (2001) 119--204.

\bibitem{Sal}
\sc D… Salamon,
\rm The Kolmogorov-Arnold-Moser theorem.
\rm Preprint, ETH-Zü\-rich, 1986.
\it Math. Phys. Electron. J. \bf 10 \rm (2004) paper 3, 37~p.

\bibitem{SZ}
\sc D… Salamon \& E…Zehnder,
\rm KAM theory in configuration space. 
\it Comment. Math. Helv. \bf{64} \rm(1989) 84--132.

\bibitem{Sev}
\sc M… Sevryuk,
\rm KAM-stable Hamiltonians.
\it J. Dyn. Control Syst. \bf{1} \rm(1995) 351--366.

\bibitem{SM}
\sc C…L… Siegel \& J…K… Moser,
\it Lectures on Celestial Mechanics.
\rm Springer, 1971.

\bibitem{St}
\sc E… Stein,
\it Singular Integrals and Differentiability Properties of
Functions. 
\rm Princeton, 1970.

\bibitem{Wh}
\sc H… Whitney,
\rm Analytic extensions of differentiable functions defined
in closed sets. 
\it Trans. A.M.S.  \bf{36} \rm(1934) 63--89.

\bibitem{XYQ}
\sc J… Xu, J… You \& Q… Qiu,
\rm Invariant tori for nearly integrable Hamiltonian systems 
with degeneracy.
\it Math. Z. \bf{226} \rm(1997) 375--387.

\bibitem{Ze73}
\sc E…Zehnder,
\rm Homoclinic points near elliptic fixed points. 
\it Comm. Pure Appl. Math. \bf{26} \rm(1973) 131--182.

\bibitem{Ze76}
\sc E…Zehnder,
\rm Generalized implicit function theorems with applications
to some small divisor problems I.
\it Comm. Pure Appl. Math. \bf{28} \rm(1975) 91--140. 
\rm II. \bf 29 \rm (1976) 49--111.

\end{thebibliography}
\end{document}